\documentclass[12pt,fleqn]{article}
\usepackage{graphicx}

\begin{document}

\newcommand{\rf}[1]{(\ref{#1})}
\newcommand{\rff}[2]{(\ref{#1}\ref{#2})}

\newcommand{\ba}{\begin{array}}
\newcommand{\ea}{\end{array}}

\newcommand{\be}{\begin{equation}}
\newcommand{\ee}{\end{equation}}

\newcommand{\const}{{\rm const}}
\newcommand{\ep}{\varepsilon}
\newcommand{\Cl}{{\cal C}}
\newcommand{\rr}{\vec r}

\newcommand{\e}{{\bf e}}

\newcommand{\m}{\left( \ba{c}}
\newcommand{\ema}{\ea \right)}
\newcommand{\mm}{\left( \ba{cc}}
\newcommand{\miv}{\left( \ba{cccc}}

\newcommand{\scal}[2]{\mbox{$\langle #1 \! \mid #2 \rangle $}} 
\newcommand{\ods}{\par \vspace{0.5cm} \par}

\newtheorem{prop}{Proposition}
\newtheorem{Th}{Theorem}  
\newtheorem{lem}{Lemma}
\newtheorem{rem}{Remark}
\newtheorem{cor}{Corollary}
\newtheorem{Def}{Definition}
\newtheorem{open}{Open problem}
\newtheorem{ex}{Example}
\newtheorem{exer}{Exercise}

\title{\bf An orbit-preserving discretization of the classical Kepler problem }

\author{
 {\bf Jan L.\ Cie\'sli\'nski}\thanks{\footnotesize 
 e-mail: \tt janek\,@\,alpha.uwb.edu.pl}
\\ {\footnotesize Uniwersytet w Bia\l ymstoku,  
Instytut Fizyki Teoretycznej, }
\\ {\footnotesize ul.\ Lipowa 41, 15-424  
Bia\l ystok, Poland}
}

\date{}

\maketitle

\begin{abstract} 
We present a remarkable discretization of the classical Kepler problem which preserves its trajectories and all integrals of motion. The points of any discrete orbit belong to an appropriate continuous trajectory. 
\end{abstract}
\ods
\ods

\ods

In last years there is a growing interest in the subject of the geometric integration of ordinary differential equations \cite{HLW}.
Geometric integration consists in numerical  solution of differential equations while preserving some physical or mathematical properties of the system exactly  (i.e., up to round-off error) \cite{LQ}. In this Letter we   propose a modification of this technique focusing not only on properties of the continuous system but also on methods of generating its exact solutions. Similar approach is well known in the case of integrable systems of nonlinear partial differential equations \cite{HA,AHS}.

In a recent paper \cite{CR} we discuss discretizations of the classical harmonic oscillator equation $\ddot x = - x$. One of these discretizations, namely
\be  \label{exos}
\frac{x_{n+1} - 2 x_n + x_{n-1}}{4 \sin^2 (h/2) }  = - x_n 
\ee
turns out to be {\it exact}. Indeed, the general solution of \rf{exos}
\be  \label{exsol}
  x_n = x_0 \cos (n h) + \frac{x_1 - x_0 \cos h}{\sin h} \sin (n h) 
\ee
satisfies  $x_n = x (h n)$, where the constant $h$ is the time step and $x (t)$ is the solution of the equation $\ddot x = - x$  with the following  initial data: $x (0) = x_0$, $\dot x (0) = (x_1 - x_0 \cos h)/\sin h $. 
It means that the points of the discrete solution $x_n$ coincide with the continuous solution $x (t)$ evaluated at the discrete time-lattice, and this equality is strictly satisfied for any $n$. 

The existence of a discretization of this kind seems to be  an exceptional phenomenon. The list of differential equations for which exact discretizations were found is rather short \cite{Ag}, although includes all ordinary differential equations with constant coefficients \cite{Ag,Reid}.
Actually, it is not easy even to find  discretizations preserving the energy integral for one-dimensional Newton equations, however in this case the problem of finding equations admitting ``integrable'' discretizations is to some extent solved  \cite{Suris}. The algorithms preserving both the symplectic form and the energy integral are often characterized by variable time stepping \cite{Lee2,KMO}.

In this Letter we consider the classical Kepler problem 
\be  \label{Kepler}
\frac{d  \vec p}{d t} = - \frac{k \rr}{r^3}  \ , \qquad 
\vec p = m \frac{d \rr}{d t} \ .
\ee
We recall that the angular momentum $\vec L$, the total energy $E$ and the Runge-Lenz vector $\vec A$ (pointing at the perihelium), namely
\be  \label{LEA}
\vec L = \rr \times \vec p \ , \quad 
E =  \frac{(\vec p)^2}{2 m}  - \frac{k}{r}  \ , \quad 
\vec A = \frac{\vec p \times \vec L}{m} - k \frac{\rr}{r} \ ,
\ee
are integrals of motion for \rf{Kepler}. The exact solution  $\vec r (t)$ is not known (except some particular cases, e.g., circular orbits) but all trajectories (orbits) can be exactly found. In polar coordinates all orbits are  described by the formula 
\be  \label{rfi}
  r = \frac{p}{1 + e \cos(\varphi - \varphi_0)} \ , 
\qquad p = \frac{L^2}{km} \ , \quad e = \sqrt{1 + \frac{2 E L^2}{m k^2}} \ ,
\ee
where $\varphi_0$ is an arbitrary constant (the initial condition).

We are going to solve the following problem: {\it find a discretization of \rf{Kepler} which preserves all integrals of motion and, if possible, preserves the trajectories}.

The existence of such discretization is far from being obvious. On the contrary, in numerical approaches to the Kepler problem it was difficult  to preserve all the integrals of motion, not saying about the orbits \cite{Le}. 
Actually, it is not easy to preserve even the energy integral \cite{Le,Is}.

Fortunatelly, a discretization of this kind exists. What is more, this algorithm is  surprisingly simple: 
\be \label{Kep-dis}
\frac{\Delta \vec p_{n}}{\Delta t_{n}}  = - \frac{ k {\vec r}_{n+1} }{\alpha r_{n+1}^2 r_{n} \cos\delta } \ , \qquad   \vec p_n := m \frac{ \Delta \rr_n}{\Delta t_n}  \ ,
\ee
where 
\[
\Delta \rr_n = \rr_{n+1} - \rr_{n} \ , \quad \Delta \vec p_n = \vec p_{n+1} - \vec p_{n} \ , \quad \Delta t_n = t_{n+1} - t_{n} \ , 
\]
and $\alpha$ and $\delta$ are constant parameters. 
The discrete time lattice $t_n$ is chosen in such a way that the angle $\Delta_n$ between $\rr_{n+1}$ and $\rr_n$ does not depend on $n$. We denote the half of this angle by $\delta_n$. Therefore, by assumption, $\delta_n = \delta = \const$.  Note that $\Delta_n = \Delta = 2 \delta$ and 
$\rr_n \cdot \rr_{n+1} = r_n r_{n+1} \cos\Delta$. 
As initial data for \rf{Kep-dis} we take: $\rr_0$, $\rr_1$  and $\Delta t_0$. The time step $\Delta t_n$ can be computed using the condition $\delta_n = \const$, compare \rf{Delta-iter}. 
In the continuum limit $\delta \rightarrow 0$. Therefore, comparing \rf{Kepler} and \rf{Kep-dis}, we see that in this limit $\alpha \rightarrow 1$.

The discrete equation \rf{Kep-dis} has several remarkable properties. They are very unusual as far as discrete analogues of continuous equations are concerned. This discretization preserves all important properties of the Kepler problem \rf{Kepler}, including its integral of motions and its trajectories: elliptic, hyperbolic and parabolic orbits. 
The presence of the parameter $\alpha$, not spoiling advantages of this discretization, can be used to obtain a better simulation of the continuous solutions by the discrete solutions. In this way we can fit, for instance, the period of the elliptical orbit.

The integrals of motion for the discrete Kepler problem \rf{Kep-dis} are given by:
\[
 \vec  L_n = \alpha \vec R_n \times \vec p_n  \ , \quad 
 E_n = \frac{(\vec p_n)^2}{2 m} - \frac{k}{\alpha R_n}  
\ , \quad
\vec A_n = \frac{\vec p_n \times \vec L_n }{m} - \frac{k \vec R_n}{R_n} \ .
\]
where 
\be  \label{R}
\vec R_n := \frac{r_{n+1} \rr_n + r_n \rr_{n+1}}{r_n + r_{n+1}} \ , \qquad   R_n := |\vec R_n | =  \frac{2 r_n r_{n+1} \cos\delta}{r_n + r_{n+1}} \ .
\ee
The vector $\vec R_n$ has a simple geometric interpretation: 
the end of $\vec R_n$ lies in the center between the ends of $\rr_n$ and $\rr_{n+1}$ (in particular, $\vec R_n$  bisects the angle between $\rr_n$ and $\rr_{n+1}$). The angular momentum can also be expressed directly in terms of $\rr_n$ and $\rr_{n+1}$
\be  \label{Lr}
 \vec L_n = \frac{\alpha m \rr_n \times \rr_{n+1}}{\Delta t_n} \ , \qquad 
L_n  \Delta t_n = 2 m \alpha r_{n+1} r_n \sin\delta \cos\delta  \ .
\ee
One can show by straightforward calculation that 
$\vec L_n$, $E_n$ and $\vec A_n$ do not depend on $n$.  Indeed, the discrete Kepler equation \rf{Kep-dis} is equivalent to
\be  \label{Kep-dis1}
\frac{\rr_{n+1}}{\Delta t_{n}}  + \frac{\rr_{n-1}}{\Delta t_{n-1}} = \left( \frac{1}{\Delta t_{n}} + \frac{1}{\Delta t_{n-1}} - \frac{k \Delta t_{n-1} }{\alpha m r_{n}^2 r_{n-1} \cos\delta } \right) \rr_{n} \ .
\ee
Therefore, 
\[
\vec L_{n} - \vec L_{n-1} = \alpha m \rr_{n} \times \left( \frac{\rr_{n+1}}{\Delta t_{n}}  + \frac{\rr_{n-1}}{\Delta t_{n-1}} \right) = 0 \ .
\]
Applying the second equation of \rf{Lr} to $L_{n} = L_{n-1}$ we get the  identity
\be  \label{ident}
 r_{n+1} \Delta t_{n-1} = r_{n-1} \Delta t_n \ .
\ee
Then, using \rf{Kep-dis} and \rf{ident}, we compute $E_n - E_{n-1} = 0$ and $\vec A_{n} - \vec A_{n-1} = 0$. All  integrals of motion can be expressed by initial data: $\rr_0$, $\rr_1$ and $\Delta t_0$. 

Now, we proceed to showing that the equation \rf{Kep-dis1} yields an explicit numerical scheme  to produce $\rr_n$ from the initial data $\rr_0, \rr_1$ and $\Delta t_0$. The crucial point is to find an appropriate iterative procedure for $\Delta t_n$. The equation \rf{ident} is not sufficient because it contains $r_{n+1}$. We have to use the constancy of $\delta_n$. Namely, $\delta_n = \delta$ implies
(after elementary geometric considerations)
\be  \label{vers}
  \frac{\rr_{n+1}}{r_{n+1}} +  \frac{\rr_{n-1}}{r_{n-1}} = 2 \cos\delta \ \frac{\rr_n}{r_n} \ .
\ee
We substitute \rf{vers} into \rf{Kep-dis1}, eliminate $r_{n+1}$ using  \rf{ident}, and as a result we obtain 
\[
\left( \frac{2 r_{n-1} \cos\delta}{\Delta t_{n-1}} - \frac{r_n}{\Delta t_n}  - \frac{r_n}{\Delta t_{n-1}} + \frac{k \Delta t_{n-1}}{m r_n r_{n-1} \alpha \cos\delta} \right) \frac{\rr_n}{r_n} = 0 \ ,
\]
and it yields 
\be  \label{Delta-iter}
\Delta t_n = \frac{\Delta t_{n-1}}{\displaystyle 2 \cos\delta \ \frac{r_{n-1}}{r_n} - 1 + \frac{k r_{n-1} (\Delta t_0)^2}{m r_1^2 r_0^2 \alpha \cos\delta}} \ ,
\ee
where we took into account $r_n r_{n-1} \Delta t_0 = r_1 r_0 \Delta t_{n-1}$ which follows from the angular momentum conservation law $L_{n-1} = L_0 $. Therefore, the complete explicit procedure of iterative computation of $\rr_{n}$ consists of two discrete equations \rf{Kep-dis1}, \rf{Delta-iter} and initial data $\rr_0$, $\rr_1$, $\Delta t_0$.

In order to obtain trajectories for the continuous Kepler problem one has to express the first two conservation laws \rf{LEA} in polar coordinates. The same approach is effective also in the discrete case. 
First, using \rf{Lr} to eliminate $\Delta t_n$, we transform the expression for the kinetic energy obtaining 
\[
\frac{(\vec p_n)^2}{2 m} = \frac{{\cal L}^2}{m^2 \alpha^2 (r_{n+1} r_n)^2} \left( \frac{r_{n+1} - r_n}{2 \sin \delta \cos\delta} \right)^2 + \frac{ {\cal L}^2}{m^2 \alpha^2 r_n r_{n+1} \cos^2\delta } \ ,
\] 
where ${\cal L} = | \vec L_n |$. 
Then, still in a direct analogy with the continuous case, 
we introduce a new variable $u_n$:
\be  \label{undef}
\displaystyle u_n := \frac{1}{r_n} - \frac{\alpha k m \cos\delta}{{\cal L}^2} \ ,
\ee
and the energy conservation law can be rewritten as 
\be  \label{osc-en-u}
\frac{(u_{n+1} - u_n)^2 }{( 2 \sin\delta)^2 } +  u_{n+1} u_n  = \frac{ 2 m {\cal E} \alpha^2 \cos^2\delta}{ {\cal L}^2 } + \frac{k^2 \alpha^2 m^2 \cos^2 \delta}{ {\cal L}^4} 
\ ,
\ee 
where ${\cal E} = E_n$ and the right hand side of \rf{osc-en-u}, although complicated, is constant (does not depend on $n$). Sometimes it is convenient to use another, equivalent, form of \rf{osc-en-u}:
\be   \label{osc-en-u1}
\left( \frac{u_{n+1} - u_n \cos\Delta}{\sin\Delta} \right)^2 + u_n^2 =  \frac{m^2 k^2 \alpha^2}{{\cal L}^4} \left( 1 + \frac{2 {\cal E} {\cal L}^2}{m k^2} \right) \ .
\ee

Both \rf{osc-en-u} and \rf{osc-en-u1} can be interpreted as the energy integral for harmonic oscillator equation. Indeed, \rf{osc-en-u} implies   
\be
\frac{u_{n+1} - 2 u_n + u_{n-1}}{(2\sin\delta)^2} + u_n = 0 
\ ,
\ee
which  can be easily solved (compare \rf{exos}, \rf{exsol}):
\be
u_n = u_0 \cos n \Delta + \frac{u_1 - u_0 \cos\Delta}{\sin\Delta} \ \sin n \Delta \ ,
\ee
and, after an elementary trigonometric transformation, we get
\be  \label{unsol}
u_n = \frac{\ep}{\cal P} \cos ( n \Delta - \theta_0)
\ee
where
\[
\tan\theta_0 = \frac{u_1 - u_0 \cos\Delta}{u_0 \sin\Delta} \ ,
\quad 
 {\cal P} = \frac{{\cal L}^2}{km \alpha} \ , \quad 
\ep = 
 \sqrt{1 + \frac{ 2 {\cal E} {\cal L}^2}{m k^2}} \ .
\]
From \rf{undef} and \rf{unsol} we obtain general expression for the orbits in the discrete case:   
\be  \label{rfi-dis}
r_n = \frac{\cal P}{\cos\delta +   \ep \cos (n \Delta - \theta_0) } 
\ . 
\ee  
This formula is similar to \rf{rfi}. However, in order to obtain the exact discretization of the trajectory we have to  require: $ {\cal P} = p \cos\delta$ and $\ep  = e \cos\delta$, which implies
\be
{\cal L} = L \sqrt{\alpha \cos\delta} \ , \quad {\cal E} = \frac{E \cos\delta}{\alpha} - \frac{m k^2 \sin^2 \delta }{2 \alpha L^2 \cos\delta} \ .
\ee

In the case of a periodic motion (negative $E$) the parameter  $\alpha$ can be fixed by the requirement that the period of the motion in the discrete case is the same as in the continuous case. 

After this work was completed, I become aware of a sequence of papers by Minesaki and Nakamura on discretizations preserving integrals of  motion \cite{MN-Kep1,MN-Kep2,IM-2cen,MN-Stack}. In particular, they also succeeded to preserve all integrals of motion and the orbits in the case of the Kepler problem. However, the method used in my paper is different from their  approach. It would be interesting to compare results generated by both discretizations.

{\it Acknowledgments.} The work was partially supported by the KBN grant No.\  1 P03B 017 28. The problem of finding the best discretization of the Kepler problem aroused in the framework of the cooperation with Bogus\l aw Ratkiewicz.


\end{document}